\input amstex
\input amsppt.sty
\magnification=1200

\parindent 20 pt
\NoBlackBoxes
\vsize=7.50in

\define \a{\alpha}

\define \dl{\delta}

\define \G{\Gamma}
\define \lm{\lambda}
\define \k{\kappa}

\define \om{\omega}
\define \r{\rho}
\define \s{\sigma}

\define \ve{\varepsilon}
\define \vp{\varphi}

\def\cupl{\operatornamewithlimits{\bigcup}\limits}

\define \fc{\frac}
\define \iy{\infty}

\define \p{\partial}

\define \ov{\overline}

\define \edm{\enddemo}
\define \ep{\endproclaim}

\define \bk{\bigskip}

\define \1{^{-1}}
\define \2{^{-2}}

\define \BC{\Bbb C}

\define \BF{\Bbb F}

\define \BN{\Bbb N}
\define \BP{\Bbb P}

\define \BZ{\Bbb Z}

\define \CA{\Cal A}

\define \CH{\Cal H}

\define \CO{\Cal O}

\define \const{\operatorname{const}}

\topmatter
 \title
 {  Affine surfaces with $AK(S)=\BC$}
 \endtitle
 \leftheadtext{T. Bandman and L. Makar-Limanov}
 \rightheadtext{ Surfaces with $C^+$ actions}
 \author
 T. Bandman and L. Makar-Limanov
 \endauthor
 \thanks The first author  is supported by the
Excellency Center of Academia and  by the Ministry of Absorption,
 State of Israel,
by  the Emmy N\"other Institute for
 Mathematics of Bar-Ilan University.
 The second  author is supported by an NSF grant.
 \endthanks

\baselineskip 20pt
\abstract
In this paper we give a  description of hypersurfaces with $AK(S)=\BC.$
Let $X$ be an affine variety and let $G(X)$ be the group generated by  
all $\BC^+$-actions on $X$. Then $AK(X)\subset \CO(X)$
is the subring  of all regular  $G(X)-$ invariant functions on $X.$ 
We give here a description of affine surfaces $S$ with $AK(S)\cong \BC$.
We show that if $AK(S)=\BC$ then $S$ is quasihomogeneous and so may be 
obtained from a smooth rational projective surface by deleting a divisor
of special form, which is called a ``zigzag''.
We denote by $\Cal A$ the set of all such surfaces, and by $\Cal H$ those
which have only three components in the zigzag. We prove that for a  
surface $S$ with $AK(S)\cong \BC$ the following statements are equivalent:
 1. $S$ is isomorphic to a hypersurface;
  2. $S$ is isomorphic to a hypersurface $S'=\{(x,y,z)\in 
\BC^3 |xy=p(z)\},$ where $p$ is a polynomial with simple roots only;
 3. $S$ admits a fixed-point free $\BC^+$- action; 
 4. $S\in\CH.$ Moreover, if $S_1\in \CH$ and
$S_2\in \CA\setminus \CH,$
then $S_1\times \BC^k\not\simeq S_2\times 
\BC^k$ for any $k\in\BN$.

\endabstract
\endtopmatter

\subheading{1. Introduction}

In this paper we  proceed with our research of the
smooth surfaces with $\BC^+$-actions  (\cite{BML1}, \cite{BML2}). 
We denote by $\CO(S)$ the ring of all regular functions on $S$. 
Let us recall that $AK$ invariant $AK(S)\subset \CO(S)$ of a surface
$S$ is just the subring of the ring $\Cal {O}(S),$ consisting of 
those regular functions on $S,$ which are
invariant under all $\BC^+$-actions of $S$. This invariant can be
also described as the subring of $\Cal {O}(S)$ of
all functions which are constants for all locally nilpotent
derivations of $\Cal {O}(S)$ (\cite{KKMLR}, \cite{KML}, \cite{ML1}).

We would like to give the answer to the following
question. 

What are the surfaces with the trivial invariant $AK$?

It is quite easy to show that the complex line $\BC$
is the only curve with the trivial invariant (\cite{ML2}).  
It is also well known, that
if $AK(S) = \BC$ and  $\CO(S)$ is a unique factorization domain,
  then $S$ is an affine complex plane 
$\BC^2$ (\cite{MiSu}, \cite{Su}). If we drop the UFD condition,
 then we have a lot of
smooth surfaces with trivial invariant, e. g. any hypersurface
of the form  $\{xy = p(z)\}\subset \BC^3,$ where all roots of $p(z)$ are
simple. 

Since we did not know any other examples, we had the following
working conjecture:

\proclaim{Conjecture}
Any smooth affine surface $S$ with $AK(S)=\BC$ is isomorphic to a hypersurface
$$\{xy=p(z)\}\subset \BC^3.$$\ep

It turned out that this conjecture is true only with an additional
assumption that $S$ admits a fixed point free $\BC^+$-action. 
Also, if we assume that $S$ is a hypersurface with $AK(S) = \BC,$
then it is indeed isomorphic to a one, defined by the equation $xy = p(z)$. 

The surfaces of this kind are well-known from the 1989 because
of the following remarkable fact which was discovered by
Danielewski (\cite {D})
in connection with the generalized Zariski
conjecture (see also Fieseler, \cite{Fi}). 
It is known that  the surfaces $\{x^ny = p(z)\}$ with $n >
1$ are not isomorphic to $\{xy = p(z)\}$
(actually they are pairwise non-isomorphic).
Nevertheless, the cylinders over all these surfaces are isomorphic. 
($S\times \BC^n$ is called ``the cylinder over surface $S$".) 
So it seems natural to
introduce a notion of equivalence for the surfaces, where two
surfaces are equivalent when cylinders over these surfaces
are isomorphic. That is why we also try to consider surfaces
with $AK (S)= \BC$ up to this equivalence. Though we are far from
complete understanding, we know that there are two classes of
surfaces which cannot be mixed by this equivalence relation.
The first class consists of the hypersurfaces $\{xy = p(z)\}$
which were mentioned above. Here is an example of a surface from
the second class:

$$S=\left\{ \matrix\format \l&&\quad\l\\  && xy=(z^2-1)z\\
(x,y,z,u)\in\BC^4: && zu=(y^2-1)y\\
 && xu=(z^2-1)(y^2-1)\endmatrix\right\}.$$
\bigskip

\subheading {2. Definitions and related notions} 

If $AK(S) = \BC$, then the group
of automorphisms of $S$ has a dense orbit. So it is natural to
compare these surfaces with quasihomogeneous surfaces which were
investigated by  M. Gizatullin, V. Danilov, and J. Bertin
(\cite{G1}, \cite{G2}, \cite{GD},
\cite{Ber}).

\definition {Definition} A smooth affine surface
 $S$ is called quasihomogeneous if the group $Aut (S)$ of
all automorphisms of $S$ has an orbit $U=S\setminus N$,
where $N$ is a finite set.\enddefinition

We will show that if $AK(S) = \BC$, then indeed $S$ 
is a  quasihomogeneous surface. Therefore,
  $S$ may be obtained from a smooth rational
 projective surface $\ov S$
by deleting a divisor
of special form, which is called a ``zigzag'' 
( \cite{G1}, \cite{G2}, \cite {GD},
\cite{Ber}).

Let us  denote by $\Cal A$ the set of all  surfaces
$S$ with $AK(S) = \BC$,
 and by $\Cal H$ those
which have only three
components in the zigzag.

We prove in Section 3, that a surface $S\in\Cal A$ is isomorphic to a
hypersurface if and only if
$S\in \Cal H$ (Theorem 1).
 In Section 4 we use this fact to prove that:

1) if
$S_1\in \Cal H$ and
$S_2\in A\setminus \Cal H,$ then the cylinders $S_1\times \BC^k$ and
$S_2\times\BC^k$ cannot be
isomorphic (Theorem 2);

2) a surface $S\in \Cal A$ admits a fixed-point
free $\BC^+$-action if
and only if $S\in \Cal H$ (Theorem 3).

Further on  the following notations will be used:

$\CO(X)$ -- the ring of regular functions on a variety $X;$ 

$K(S)$ -- canonical divisor of a surface $S;$

$[D]$ --  class of linear equivalence of a divisor $D$;

$\tilde D$ -- proper transform of a divisor $D$ after a  blow-up;

$ D^*$ --   algebraic (total) transform of a divisor $D$ after 
a blow-up;

$(\omega), \ (f)$ --  divisors of zeros of a form $\omega$ and a function
$f$, respectively;

$Aut(S)$ --  automorphism group of a surface $S$;

$G(S)$ -- subgroup of $Aut(S),$ generated by all $C^+$-actions 
on a surface $S$;

$OG(S)$ -- a general orbit of the group $G(S)$.

$\ov A$ - a  Zariski closure of $A$ if other meaning is not specified.

`` General'' means `` belonging to a Zariski open subset''

A singular point of a rational function is a point where 
 the function  is not defined.
\bk

\subheading{3.\ Characterization of hypersurfaces $S$ with $AK(S)=\BC$}

Following \cite{Ber}, \cite{Mi1}, \cite{MiSu}, by a line pencil 
on a surface $S$ we mean a morphism $\rho: S\to C$  into a smooth curve $C,$
such that the fiber  $\rho^{-1} (z) $ for a general $z\in \BC$ is 
isomorphic to $\BC.$ 
 Then $S$ contains a cylinder-like subset i. e. an open subset which
is isomorphic to a direct product of $\BC$ and an open subset of $C$ 
(\cite{B}, III.4).
The pencils are different if their general fibers do not coincide.  
Any line pencil $\rho$  over  affine curve  $C$
on a surface $S$ corresponds anto a $C^+$-action
$\vp_{\rho}$ on $S,$ such that its general orbit coincides with a 
general fiber
of the pencil, and  corresponds to a locally nilpotent derivation ($lnd$)
$\partial _\rho$ in
the ring $O(S)$ of regular functions on $S,$ such that $\partial _\rho f=0$
if and only if $f$
is $\vp_\rho$-invariant (\cite{ML1},\cite{KML}, \cite{S}, \cite{Mi1}).
If there are two different  line pencils in $S,$ then $\rho(S)=\BC$ 
  (indeed, in this case $\rho(S)$
is an affine curve containing the image of a fiber of 
the second line pencil, and this fiber is 
isomorphic to $\BC$). 
 Since we are looking for the surfaces having many $\BC^+$ actions, 
we assume further on that $C\cong\BC$.

For a   pencil $\rho$ over $\BC$
one can find such a closure $\ov S$ of $S$ that the
extension $\ov\rho:\ov
S\to\Bbb P^1$ of the map $\rho: S\to\BC$ is regular and in the
commutative diagram
$$\alignat3
&\quad S\quad &&\hookrightarrow\quad &&\ov S\\
&\rho \downarrow && && \downarrow {\ov\rho}\tag1\\
&\quad \BC\quad &&\hookrightarrow\quad &&\Bbb P^1\endalignat$$
the divisor $B=\ov S\setminus S$ is connected and  has the following properties:

\noindent I.\quad $B=F+D+E,$ where

a)\ $F\cong \Bbb P^1,$ \ $\ov\rho(F)=\BP^1-\BC,$

b)\ $\ov\rho\bigm|_D:D\to\BP^1$ is an isomorphism,

c)\ $E=\sum E_i+\sum H_i,$ where
$\ov\rho(H_i)\in \BC\setminus \rho(S)$
and
 $\ov\rho(E_i)=z_i\in \r(S)$ is a point.
Moreover, $\ov{\rho^{-1}(z_i)}$ is a union of disjoint 
smooth rational curves, and each of
them intersects $B$ precisely at one point.

\noindent II.\quad $B$ does not contain $(-1)$ curves, except, maybe, $D$.

 The
structure of fibers is described in \cite{Mi1}, Lemma 4.4.1. If
 there are two different  line pencils in $S,$ then 
 $E=\sum E_i.$

\definition{Definition}
We call a closure $\ov S$ a good $\rho$-closure of an affine surface $S,$
if it has properties
I and II.\enddefinition

\definition {Definition}  Let $F_z=\r^{-1}(z)=\sum _{i=1}^{i=m}n_iC_i$,
 where $C_i$ are connected (and irreducible) components.
 If  $m=1, n_1=1,$ then the fiber is called non-singular. 
The singular fiber  is either non-connected or has $m=1,n_1>1.$
\enddefinition

\bigskip
\proclaim{ Proposition 1} Let $S$ be a smooth affine surface
  with a line pencil $\rho.$ 
  Let $\ov S$ be
a good  $\rho-$closure of $S.$ Let $F_{z_1},..., F_{z_n}$
be all singular fibers of $\rho,$ and let 
$F_{z_i}=\sum _{j=1}^{j=k_{i}} n_{i,j}C_{i,j}$ 
be a sum of irreducible curves $C_{i,j}; \ \ C_{i,j}\cong \BC.$
 Then there exists a function $\a \in \CO ( S)$, such that 

a) $\a$ is linear along each nonsingular fiber $F_{z}, z\ne z_i, i=1,...,n$ 
(i.e. $\a\bigm |_{F_{z}}$ is  a non-constant linear function);

b) $\a \bigm |_{C_{i,j}}=\a_{i,j}=const$ for all $1\leq i\leq n,
1\leq j\leq k_i$.\ep

\demo{Proof of  Proposition 1}
Let $\partial_{\rho}$ be a non-zero $lnd,$
 corresponding to the line pencil $\rho.$
If there is a nonsingular fiber $F_z=\r^{-1}(z),$ such that
 $\p_\r(v)\bigm|_{F_{z}}=0$ for all $v\in\CO(S),$ then we may 
consider another $lnd$ \ $\tilde\p_\r=\fc{\p_\r}{\r-z}$
and repeat  this procedure, if needed. Hence we may assume that
$\p_{\r}$
 does not vanish identically along the nonsingular fibers of $\r.$

 Since $\partial _\rho$ is a non-zero derivation,
 there exists a function $v\in\CO(S),$ for which $\partial _\rho (v)\ne 0,$
i.e. the  minimal $n$ for which $\partial _\rho^n (v)= 0$
is not smaller then $2.$ Let us take $u= \partial _\rho^{n-2} (v).$
Since $\partial _\rho^2 (u)= 0,$  \ $\partial _\rho (u)=f(z)$
    depends only on $z=\r(s), s\in S.$
If $f(\tilde z)=0, \quad \tilde z\ne z_1,...,z_n,$ then 
$u\big |_ {\r^{-1}(\tilde z)}=u_0={\const},$ and  we
consider a new
 function $\fc{u-u_0}{\r-\tilde z}.$

Repeating this, we arrive at 
the situation, such that
\roster\item $\partial_\rho u= f(z)$, where $f$ vanishes only
at the points $z_i, i=1,...,n;$
\item $u$ is a linear function along each fiber $\rho^{-1}(\tilde z),$
  \ $\tilde{z}\ne z_i, i=1,...,n.$\endroster
\bigskip

We are going to show, that
 $u=u_i=\const$ along each component  $C_{i,j}$ of
$F_{z_{i}}, \   i=1,...,n.$

Indeed, $u$ is linear along a general fiber, which means that the
intersection $(\ov{U}_w,\ov {\r}^{-1}(z))=1$
 for closure $\ov{U}_w,$
in $\ov S$ of a 
a general level curve $U_w=\{s\in S: u(s)=w\}$ and any $z.$

If $u\bigm|_{C_{i,j}}\ne const,$ then $(\ov{U}_w, C_{i,j})\geq 1,$
and  $(\ov{U}_w,\ov {\r}^{-1}(z_i))\geq n_{i,j}. $ Thus, if 
$n_{i,j}>1,$ then  $(\ov{U}_w, C_{i,j})=0$ and $u\bigm|_{C_{i,j}}= const.$

If $n_{i,j}=1$, then the fiber is non-connected
and $u\bigm|_{C_{i,j}}\ne const$ implies that $U_w$
does not intersect $\ov{\r}^{-1}(z_{i})\setminus C_{i,j}$
for a general $w\in \BC.$ Thus,
 $u \bigm|_{\ov{\r}^{-1}(z_i)\setminus C_{i,j}}$
 has to be regular and constant. On the other hand $u$ has a 
pole along $D,$ so  
$u \bigm|_{\ov{\r}^{-1}(z_{i})\setminus C_{i,j}}=\infty.$
Since $u$ has only regular points, $u\bigm|_{C_{i,k}}=\infty$
 as well, if $k\ne j$.   But
$u\in\CO(S)$ so there are no 
components with $k\ne j.$ Therefore  $\rho^{-1}(z_i)$
has just one component of multiplicity 1,
which contradicts  to our  assumption.

Thus, we may take
 $\a=u.$
\hfill$\qed$\edm

\bigskip

\proclaim{ Proposition 2}
Any smooth affine surface $S$ with $AK(S)\cong \BC$
is quasihomogeneous.\ep

\demo {Proof of Proposition 2} Assume that $\phi$ and $\psi$
are  $ C^+$-actions on $S$ having different orbits.
Let $\r$ and $\k$ be  the corresponding line pencils,
$\p_\r,\ \p_\k$ corresponding $lnd,$
 let
$R_z=\r^{-1}(z)$ and $K_w=\k^{-1}(w)$ for general $z,w \in\BC$
and let $\ov {R}_z$ and $\ov {K}_w$ be their closures in a good 
$\r-$closure  $\ov S$ of $S.$

We are going to show that $S\setminus OG(S)$ is a finite set.

If a point $s\in S\setminus OG(S),$ and the fiber $R_{\r(s)}$ 
is non-singular, than  $R_{\r(s)}\subset S \setminus OG(S)$
as well. Indeed, as it was shown in Proposition 1, 
 we can choose $\p_\r,\ \p_\k$
in such a way, that they do not vanish along non-singular fibers,
i.e. there are no fixed points in these fibers. 

For the same reason, $R_{\r(s)}$ does not intersect a general
fiber $K_w,$ i.e. it is contained in $ K_{\k(s)} .$

But then $\r\ne \r(s) $ along  a general $K_w,$ hence 
$\r \bigm|_{K_w}=const,$ and the fibers of these two 
actions coincide.

Thus, $s\in S\setminus OG(S)$ implies that $s\in R_{z_0}\cap K_{w_0}$
for singular fibers $R_{z_0},K_{w_0}.$

If $S\setminus OG(S)$ is infinite, then 
there exists
   a connected component
 $C\subset R_{z_0}\cap K_{w_0}$ for singular  
 fibers $R_{z_0}, K_{w_0}$ of $\r$ and $\k,$  respectively.

 Let $\ov{\r}^{-1}(z_0)= \ov C \cup E'\cup (\cup \ov{C}_i),$ where
 $E'\subset \ov S \setminus S$ and $C_i$ are
 others components of $\r^{-1}(z_0).$
Consider $K_w\cong \BC.$ Since  the intersection
$(\ov {K}_w, \ov{R}_z)\geq 1,$  \ $\ov {K}_w$ intersects
 $R_{\infty}=\ov{\r}^{-1}(\infty).$
Hence,  the only puncture of $K_w$ belongs to $R_{\infty},$
and that means that  $\ov {K}_w\cap E' =\emptyset.$
So,  $\k$ has no singular points and has to be constant along 
$E'.$ Since $E'\cap D \ne \emptyset, \k \bigm |_{E'}=\k \bigm|_{D}$ 
(see diagram 1 and recall that $E'$ is connected). But $\k \bigm |_{D}=\infty,$  (if it was not so, 
$\k$ would be bounded and, hence, constant along  a general fiber $R_z$).

We got that $\k \bigm |_{E'}=\infty$ and has no singular points.

On the other hand $\k $ is  finite and constant along  $C,$
which implies that the point $C\cap E'$ is singular.

 The contradiction
 shows that no such curve $C$ exists  and 
$S\setminus OG(S)$ is a finite set. So $S$ is indeed 
quasihomogeneous.  
\hfill $\qed$\edm

 \bigskip

Any good $\rho$-closure $\ov S$ of $S$ may be described by a 
graph $\G(\ov
S)$ in the following
way: the vertices of this graph are in bijection with irreducible 
components
of the divisor $\ov
B=\ov S\setminus S,$ and two vertices are connected by an edge, if they
intersect each other.

\bigskip
Now we are going to  use  the description of  quasihomogeneous 
affine surfaces  due to M.H.
Gizatullin and J. Bertin (\cite{G1}, \cite{G2}, \cite{GD}, \cite{Ber}).

Any such surface $S\ne\BC^2$  may be obtained by the following blow-up process,
described in \cite{G2}.

Let $S_0=\BP^1\times\BP^1.$
Let $\ov{\rho}:\BP^1\times\BP^1\to\BP^1$ be a
 projection onto the second factor.
Let $F_0=\ov{\rho}^{-1}(z_0),$\ $F_1=\ov{\rho}^{-1}(z_1),$\
 $z_0,z_1\in\BP^1$ and let
$D$ be a section; that
is, $\ov{\rho}\bigm|_D: D\to\BP^1$ is an isomorphism.

{\it Step 0} is an initial step, we start with the divisor,
 described by the
 following graph:

\noindent \quad $ \underset f\to{_\bullet} \underline{\qquad\quad}
\underset d\to{_\bullet}\underline{\qquad\quad}
\underset f_1\to{_\bullet},$
\flushpar
where vertices $f,d,f_1$ represent components $ F_0,
 D,
F_1$ respectively,

{\it Step 1} is the blow-up $\s_1:S_0\to\ov S_1$ of a point
$w_1 \in F_1$ into
an exceptional component $E\subset \ov S_1.$

We denote  $F_1^*= \tilde F_1+E$ as $E_0+E_1,$
where $E_0, E_1$ are two rational curves, and the graph
of $F_0+D+E_1+E_0 $ looks like 

\noindent \quad $ \underset f\to{_\bullet} \underline{\qquad\quad}
\underset d\to{_\bullet}\underline{\qquad\quad}
\underset e_1\to{_\bullet}\underline{\qquad\quad}
\underset e_0\to{_\bullet},$
\flushpar
where  the vertices $f,d,e_1,e_0,$ represent  the 
components $\tilde F_0,
\tilde D, 
E_1, E_0$ respectively,

{\it Step 2} is the blow-up $\s_2:\ov S_1\to\ov S_2$ of  a
 point $w_2\in
E_1\cup E_0$ into a component $E_2\subset \ov S_2,$
in such a way that a graph of  $\tilde{F}_0+\tilde{D}+
\tilde{E}_1+\tilde {E}_0 +E_2 $
is linear.

{\it Step $k$} is the blow-up $\s_k:\ov S_{k-1}\to\ov S_k$ of a
 point
$w_k\in \tilde {E}_0\cup \tilde {E}_1\cup
\tilde {E}_2\dots\cup \tilde {E}_{k-2}\cup {E}_{k-1}$ into a
component $E_k\subset \ov {S}_k$ in such a way
 that the graph of the divisor $\tilde {F}_0+\tilde
{D}+ \tilde {E}_0+\tilde {E}_1+\tilde {E}_2+\dots
+\tilde {E}_{k-1}+E_k$ is linear.

{\it Step $k+1$} is the last step.
Let $\a_1\dots\a_q$ be different points in $\tilde E_0\cup\tilde E_1
\cup\dots\cup\tilde
E_{k-1}\cup
E_k$ such that
 each $\a_i$ belongs to one component only, $1\le i\le q.$
Let $\tau_1\dots \tau_q$ be blow-ups of the points $\a_1\dots \a_q$
 into the curves $G_i,$ \ $1\le i\le q,$ respectively; and let
$\ov S$ be $\tau_1\circ \tau_2\dots\circ \tau_q$ $(\ov S_k).$

The desired surface $S=\ov S\setminus\left\{\tilde F_0\cup
\tilde D\cup\left(\cupl_{j=0}^k\tilde
E_j\right)\right\}.$

\proclaim{Remark} This description of quasihomogeneous 
surfaces implies, in particular, that there may be only one 
singular fiber  for a line pencil $\r.$\ep

We want to choose the ``minimal" way  to obtain $S$ by the 
described process,
i.e. to obtain a good $\r-$closure of $S.$
 For this we want to substitute $ S_0=\BP^1\times\BP^1$ 
by  a minimal ruled surface $\BF_n$ (see \cite{B}).

In the sequel, for simplicity of notation we will denote $\tilde E_j$ as
$E_j,$ since it cannot lead to confusion.

\proclaim{Proposition 3} The surface $S\not\cong\BC^2$, obtained 
by the blow-up process, described above, may be obtained by 
the same process, starting with the minimal surface $S_0=\BF_n$ 
and ending with such $\ov S,$
 that $E_j^2\ne {-1}$    in $\ov S$ for all $0\le j\le k.$\ep

\demo {\it Proof of Proposition 3}  We prove the Proposition by
 induction on  the number of steps $k.$ We start with the surface
 $S_0=\BF_n$ and show, that, changing $n,$ we may always eliminate 
the $(-1)$ components.

Assume that $k=0.$ Since $\rho^{-1}(z_1)\subset S$ is  singular
(recall that $S\not\cong\BC^2$)
there are  points 
$\a_i\in F_1, \ 1\le i\le q,$ in $F_1,$ which are blown up 
at the first (and last) step into the curves $G_i.$ Thus,
in $\ov S$ this fiber has a form $\tilde {F}_1+\sum_{i=1}^{i=q} G_i$
( the multiplicities are equal to $1$), which implies that the fiber
is not connected, $q>1, (\tilde {F}_1)^2=-q<-1.$

Assume now that the Proposition is true for all $k<k_0.$
Let $E_j$ be  a component of $F_1^*$ in $\ov{S}_{k_0},$
such that $E_j^2=-1.$

  There are two possibilities:

1.  $E_j$ is a result of the  blow up $\s_j.$
The 
points of this component are  not blown-up at any proceeding step,
 since it would make $E_j^2<-1.$ Thus, $E_j$ may be contracted back,
 and we may obtain surface
$S$ by the same process, omitting the step number $j,$ i.e.
 as a complement to zigzag,
 obtained
 by the blow-up process having one step less.

2. $E_j$ is a proper transform of $F_1.$  In this case we may
 blow it down after step 1, and obtain the same surface by 
the same process, having one step less,  starting with the surface
 $S_0=\BF_{n+1}$ or   $S_0=\BF_{n-1}.$

By the assumption of the induction it means that the Proposition
 is true for $k_0.$
\hfill$\qed$\edm

\definition{Definition} We denote by $\CA$ the class of all smooth 
affine surfaces $S$ with $AK(S)=\BC.$ Let us denote by 
$\CH\subset\CA$ the
subset of those surfaces, for which $k=0$
in a good $\r-$closure, obtained by the described process.
\enddefinition

\proclaim{Theorem 1}
A surface $S\in \CA$ is isomorphic to a hypersurface
 if and only if $S\in \CH.$
\ep

\demo{Proof of Theorem 1}
The proof is based on a property of hypersurfaces, which was explained to
the authors by V. Lin
and M. Zaidenberg.
We give here a proof because of lack of a reference.

\proclaim{Lemma 1}
Let $X\subset\BC^n, \, n>2$ be a smooth  hypersurface.
Then the canonical class $K(X)$ of $X$ is trivial (i.e., the divisor of
zeros of a holomorphic
$(n-1)$-form on $X$ is equivalent to zero).\ep

\demo{Proof of Lemma 1}
Let $\{z_1,\dots, z_n\}$ be coordinates in $\BC^n$ and let $p(z_1\dots
z_n)=0$ be the
equation of $X.$
Since $X$ is  irreducible and smooth,  
$p(z_1\dots z_n)$ is an irreducible polynomial 
 and the 
vector $\triangledown p=\left(\fc{\partial p}{\partial
z_1},\dots,\fc{\partial p}{\partial z_n}\right)$  does not vanish
on $X$.

Let $\omega=\sum\limits_{i=1}^n(-1)^{n+i}dz_1\wedge\overset{\overset i\to
\vee}\to \dots \wedge d
z_n,$ where $\overset i\to\vee$ denotes that $dz_i$ is omitted.
Let $\eta=\omega\bigm |_X.$
We want to show that the divisor  $(\eta)$ of the 
form $\eta$ is
linearly equivalent
to zero.

Choose a point $x\in X$ and a linear change of coordinates
$(u_1,\dots,u_n),$ such that
$u_1,\dots, u_{n-1}$ are coordinates in the tangent hyperplane to $X$
at the point $x$ and
$$u_n(z_1,\dots,z_n)=\triangledown p\pmatrix z_1\\ \vdots\\
z_n\endpmatrix=\sum_{i=1}^n\fc{\partial p}{\partial z_i}(x)\cdot z_i.$$
We denote by $A$ the matrix
$$ A=\bmatrix \frac{\partial z_1}{\partial u_1} & \cdots & \fc{\partial
z_n}{\partial
u_1}\\
\hdotsfor3\\
\fc{\partial z_1}{\partial u_n} &\cdots &\fc{\partial z_n}{\partial
u_n}\endbmatrix.$$

In local coordinates $(u_1,\dots, u_{n-1})$ in a neighborhood of $x$,
$$\eta=\sum_{i=1}^n  A_{ni}du_1\wedge\dots\wedge du_{n-1},$$
where $ A_{ni}$ are the algebraic complements of $\fc{\partial
z_i}{\partial u_n}$ in the
matrix $ A,$ i.e., $A_{ni}=\det A\cdot a_{in},$ where
$$\bmatrix a_{11} & \dots & a_{1n}\\
\hdotsfor3\\
 a_{n1} &\cdots &a_{nn}\endbmatrix =
A^{-1}=\bmatrix
\fc{\partial u_1}{\partial z_1} &\cdots& \fc{\partial u_n}{\partial z_1}\\
\hdotsfor3\\
\fc{\partial u_1}{\partial z_n} & \cdots & \fc{\partial u_n}{\partial
z_n}\endbmatrix.$$

Thus, $\eta=\left(\sum\limits_{i=1}^n\fc{\partial p}{\partial
z_i}\right)\det A\cdot
du_1\wedge\dots\wedge du_{n-1}$ and the divisor $(\eta)_0$ is the zero
divisor of a regular
function $g(x)=\sum\limits_{i=1}^n\fc{\partial p}{\partial z_i}$ on $X$.
It follows that a holomorphic $(n-1)$-form $\fc{\eta(x)}{g(x)}$ does not
vanish on $X,$ hence
$[K(X)]=0.$\hfill $\qed$\edm

Let $S\in \Cal A, S\ne \BC^2.$ 
The graph $\G(\ov S)$ has a form:

\noindent \quad $ \underset f\to{_\bullet} \underline{\qquad\quad}
\underset d\to{_\bullet}\underline{\qquad\quad}
\underset e_{t_1}\to{_\bullet}\ldots \underset e_1\to{_\bullet}\ldots \underset
e_0\to{_\bullet}\ldots\underset e_{t_k}\to{_\bullet}.$

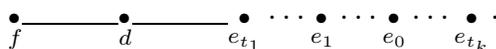
\flushpar

or the form (if $k=0$)

\noindent \quad $ \underset f\to{_\bullet} \underline{\qquad\quad}
\underset d\to{_\bullet}\underline{\qquad\quad}
\underset f_1\to{_\bullet},$

\flushpar 
where the vertices $f,d,f_1,e_1,e_0,$ represent the 
components $\tilde F_0, \tilde D, \tilde F_1, E_1, E_0$ 
respectively, and a vertex $e_{t_i}$ represents the
component $ E_{t_i}$,
obtained  at the  step $t_i$.

\proclaim{Definition}
We say that $e_i<e_j$ $(E_i<E_j)$ if $e_i$ is on the left of $e_j$ in 
the graph $\G(S).$\ep

\proclaim{Lemma 2}
The canonical  class $[K(\ov S_k)]$ of $\quad\ov S_k, \quad k>0,\quad$
 is the class of the divisor
$$K(\ov S_k)=\a\tilde F_0-2\tilde D-E_1+\sum_{i=2}^k \ve_i E_i,\tag2$$
where 

$$\a\in\BZ,\quad \ve_i<-1 \quad\text{if}\quad  
e_i<e_1, \quad\text{and}\quad
\ve_i\ge0\quad\text {if}\quad e_i>e_1.\tag 3$$
Let $$F_1^k=F_1^*=\sum_{i=0}^{i=k} n_iE_i$$ be the algebraic (total) 
transform of $F_1$ in $S_k$. Then

$$\ve_i<n_i-1 \quad\text{if}\quad e_i<e_0,\quad \ve_i\ge n_i \quad\text{if}
\quad e_i>e_0,\quad n_1=n_0=1.\tag 4$$
\ep

\demo{Proof of Lemma 2}

We prove first inequalities (3)  by induction on $k.$

The canonical class of  $\BF_n$ is $[\a F_0-2D]$
(\cite{B}, Prop.III.18).
Consider the first step: the fiber $F_1\subset \BF_n$ is blown up 
into two rational curves $F_1^*= \tilde{F}_1 +E$. 
Both curves have selfintersection $-1.$

Two cases are possible:

1. $\tilde {F}_1\cap \tilde D=\emptyset, E\cap\tilde D\neq \emptyset.$ 

According to the formula for the canonical class of a blow-up 
(\cite{Ha}, Prop. 3.3, ch. V)
 the canonical divisor

$$\align K(\ov S_1)&=\sigma_1^*(K(\BF_n))+E\\
&=\a\tilde F_0-2\tilde D-2E+E=\a\tilde{F}_0-2\tilde{D}_0-E.\endalign$$

In this case we denote $E=E_1, \tilde {F}_1=E_0.$

2. $\tilde {F}_1\cap \tilde D\ne\emptyset, E\cap\tilde D= \emptyset.$

Then the canonical divisor

$$\align K(\ov S_1)&=\sigma_1^*(K(\BF_n))+E\\
&=\a\tilde F_0-2\tilde D+E=(\a+1)\tilde F_0
-2\tilde D-\tilde{F}_1.\endalign$$

since $F_0\cong E+\tilde{F}_1.$
In this case we denote $E=E_0, \tilde {F}_1=E_1.$

Thus, for $k=1$ the formula is proved.

Assume now that (2)  and (3) are proved for all $k<k_0:$

$$K(\ov S_{k_0-1})=\a\tilde F_0-2\tilde D-E_1+
\sum_{i=2}^{k_0-1}\ve_i E_i.$$

 Then

$$\align K(\ov S_{k_0})& =\sigma_{k_0}^*(K(\ov S_{k_0-1}))+E_{k_0}\\
&=\a\tilde F_0-2\tilde D-E_1+\sum_{i=2}^{k_0-1}\ve_i\cdot
E_i+\ve_{k_0}E_{k_0}.\endalign$$

Consider the following cases:

\noindent I.\quad At the step $k_0$ we blow up a point $w_{k_0},$ belonging only
to the component
$E_s,$ represented by the vertex on the far right 
($e_s\geq e_r$
for all $r< k_0$).
In this case, $e_s$ is on the right of $e_1.$  By the induction
 assumption, we  have $\ve_s\ge 0,$ and $\ve_{k_0}=(\ve_s+1)>0.$

\noindent II.\quad At the step  $k_0$ we blow up the meeting point
$E_s\cap E_{s'},$ where
$e_s <e_{s'}\le e_1.$ 
 Then   $\ve_s<-1,\ \ve_{s'}\le-1,\ \quad \ve_{k_0}=\ve_s+\ve_{s'}+1<-1-1+1<-1.$

\noindent III.\quad At the step  $k_0$ we blow up the meeting point
$E_s\cap E_{s'}$, where
$e_s >e_{s'}\ge e_1.$ 
 Then $\ve_s\ge0,\ve_{s'}\ge-1,\ \quad \ve_{k_0}=\ve_s+\ve_{s'}+1\ge -1+1\ge 0.$

\noindent IV.\quad At the step number $k_0$ we blow up the meeting point
$E_s\cap \tilde D.$ Then
$e_s\leq e_1,\quad \ve_{k_0}=\ve_s-2+1\le -1-1< -1.$ 

Since the graph $\G(S)$ is linear, we exhausted all the possibilities.

 Now  let us prove the inequalities (4).

For $k=1$ we have $F^1_1=E_1+E_0,\quad K(\ov{S}_1)=
\a \tilde F_0 -2\tilde D -E_1$, 
therefore, $\ve_1<n_1-1.$

We prove (4) for any $k$ by induction. 
Assume that it is proved for all $k< k_0.$ 
Then in $\ov{S}_{k_0}$ we have 

$$F_1^{k_0}=\sigma_{k_0}^*(F_1^{k_0-1})=\sum_{i=0}^{i=k_0-1}n_iE_i 
+n_{k_0}E_{k_0},$$
where $n_{k_0}=n_s+n_r$, if $E_{k_0}$ appears as a blow-up of 
the intersection $E_s\cap E_r,$ 
 or $n_{k_0}=n_s$ if $E_{k_0}$ is the result
of a blow-up of  either $D\cap E_s$  or of a point of the
maximal component.

 Using the inequalities (4) for $k<k_0,$ we get:

$n_{k_0}=n_s\le\ve_s<\ve_s+1=\ve_{k_0},$ if $E_s$ is the maximal 
component  and $s \neq 0$;

 $n_{k_0}=n_0 = 1 \le 1 = \ve_{k_0}$, if $E_0$ is the maximal 
component;

$n_{k_0}=n_s+n_r\le \ve_s+\ve_r<\ve_s+\ve_r+1=\ve_{k_0},$ 
if $e_0<e_s<e_r;$

$ n_{k_0}=n_0+n_r=1+n_r\le 0+\ve_r+1=\ve_{k_0},$ if $e_0=e_s<e_r;$

$ n_{k_0}=n_s+n_0=1+n_s>1+\ve_s+1 =\ve_{k_0}+1,$ if $e_s<e_r=e_0;$

$ n_{k_0}=n_s+n_r>\ve_s+1+\ve_r+1 =\ve_{k_0}+1,$ if $e_s<e_r<e_0;$

$n_{k_0}=n_s>\ve_s+1=\ve_{k_0}+2>\ve_{k_0}+1,$ if $E_s$ is the minimal component.
\hfill $\qed$\edm

\proclaim{Lemma 3} 
Denote the transform of $F_1$ in $\ov S$
$$ F_1^{k+1}=F_1^*=\sum_{i=0}^{i=k} n_iE_i+\sum_{i=1}^{i=q}g_iG_i,$$
where $n_1=n_0=1, g_i>0, n_i>0.$

 $[K(S)]=0$ if and only if the divisor $K(\ov S)$ is
equivalent to a linear
combination 
$$\sum_{i=0}^{i=k} \a_iE_i+f\tilde {F}_0+d\tilde{D}
+m\sum_{i=1}^{i=q} g_i G_i,\tag 5 $$
for  some $m\in \BZ.$
\ep

\demo{Proof of Lemma 3}

 $$\align K(\ov S)&= K(\ov S_k)^*+\sum G_i\\
&=\a\tilde F_0-2\tilde
D-E_1+\sum_{i=1}^k\ve_iE_i+\sum_{i=1}^q\dl_iG_i,\tag 6\endalign$$

where $\dl_i=\ve_s+1$ for each $G_i$ intersecting $E_s.$

If $[K(S)]=0,$ then $K(S)$ is  divisor 
of a rational function $h,$ which has zeros and poles
in $S$ only along components $G_i.$
 But then $h$
does not vanish and has no poles in any fiber $F_z, \ z\ne z_1.$
Since  general fiber is isomorphic to $\BC,$
it means, that $h$ is constant along each fiber, 
i.e. $h(s)=(\r(s)-z_1)^m$. But then $\dl_i=mg_i.$ 
\hfill $\qed$ \edm

\proclaim{Definition} 
We call component $E_s$ essential, if there 
is a component $G_{i_s}$ of the fiber $F_1^*\subset \ov{S}$ such that 
 $G_{i_s}\cap E_s\ne
\emptyset.$\ep

\proclaim { Remark} We see from the previous Lemma that
$[K(S)]=0$ implies
 $\ve_s+1 = mn_s$
for any essential component $E_s$. \ep

\proclaim{Lemma 4}
If $k>0,$ then $[K(S)]\ne0.$\ep

\demo{Proof of Lemma 4}  Consider the graph

\noindent \quad $ \underset f\to{_\bullet} \underline{\qquad\quad}
\underset d\to{_\bullet}\underline{\qquad\quad}
\underset e_{t_1}\to{_\bullet}\ldots \underset e_1\to{_\bullet}\ldots \underset
e_0\to{_\bullet}\ldots\underset e_{t_k}\to{_\bullet}.$

 Assume that $[K(S)]=0,$ i.e
 $\ve_s+1 = mn_s.$
Several cases, concerning the place of essential components in
the graph, are possible.

I. There is an essential component $E_s$ such that $e_s \ge e_0$.
Then according to Lemma 2 $n_s \le \ve_s+1 = mn_s$ and so $m \ge 1$.

II. There is an essential component $E_s$ such that $e_1 < e_s < e_0.$ 
Then according to Lemma 2 $n_s > \ve_s+1 = mn_s > 0$ and hence 
$1 > m > 0$.

III. There is an essential component $E_s$ such that $e_s \le e_1$.
Then according to Lemma 2 $0 \ge \ve_s+1 = mn_s$ and $m \le 0$.

Therefore we may have only one of these cases for all essential 
components. 

Let us assume that $e_s \le e_1$ for any essential component $E_s.$ 
Let $t_0 = max\{t:e_t>e_1, t \ge 0\}$. By the construction 
$(E_t)^2=-1$ in $\ov{S}_k,$ (it is the result of a blow-up).
So it should contain a point
 which is blown up  at the last $k+1$ step. But then $E_t$
is essential, which is impossible in this case. 

The case $e_s \ge e_0$ for all essential components can be treated
analogously since the last component to the left of $E_0$ also 
must be essential.

The remaining case II is impossible since $m \in \BZ$.
 
Therefore (5) can be true only if the graph has three components:

\noindent \quad $ \underset f\to{_\bullet} \underline{\qquad\quad}
\underset d\to{_\bullet}\underline{\qquad\quad}
\underset f_1\to{_\bullet}.$
\hfill $\qed$\edm

\proclaim{Lemma 5} If $k=0,$ then $S$ is a hypersurface.\ep

\demo{Proof of Lemma 5}
Let $\rho: S\to\BC$ be a line pencil in $S,$\ $\ov \rho$ its extension to a
good $\rho$-closure
$\ov S$ of $S$, $\vp_\rho$ and $\partial_\rho$ the corresponding
$\BC^+$-action and $lnd$
respectively.
Let $\rho^{-1}(0)$ be the only singular fiber. 
All the multiplicities are 1 in this case, so the fiber cannot be connected.
Let $u\in O(S)$ be a function such that 
\roster\item $\partial_\rho u=\rho^n;$
\item $u$ is a linear function along each 
fiber $\rho^{-1}(z),$\ $z\ne0;$
\item $u=u_i=\const$ along each component 
 $G_i$ of $\rho^{-1}(0),$\
$i=1,\dots,q.$\endroster

Such a function exists due to Proposition 1.
We are going to show that we can choose $u,$  so that $u_i\ne u_j,$ 
when $i\ne j,$ and the rational
extension $\ov u$  of $u$ to
$\ov S$ is finite and
non-constant along $\tilde{F}_1.$ 
Indeed, $u$ is linear along a general fiber, which means that the
intersection $(\ov {U}_w,\ov{F}_z)=1$ for the closures of 
a general level curve $U_w=\{s\in S: u(s)=w\}$ and 
 the closure $\ov{F}_z$ of a
 general fiber
$F_z=\{s\in
S:\rho(s)=z\}.$

There are three possibilities:

\noindent I.\quad $\ov u\bigm|_{\tilde {F}_1}=u_0\in\BC,$ 
and $u_0\ne u_1=\ov u\bigm|
G_1.$
Then the intersection $G_1\cap \tilde {F}_1=\a_1$ is a singular point, and a general
level curve passes
through $\a_1.$
Another singular  point $\a_2=D\cap \tilde{F}_1,$ since $\ov u\bigm|_D=\iy.$
Thus, a general level curve $U_w$ has to pass through $\a_2$ as well.
But this contradicts to $(\ov{U}_w,\ov{F}_z)=1.$

Thus, $\ov u\bigm|_{\tilde {F}_1}=u_0\in\BC$ implies  
 $u_0=u_1=u_2=\cdots=u_q$ and we can consider
 a new function $\fc{u-u_0}{\rho}$ instead of $u$
(because $F_1^*=\tilde{F}_1+\sum G_i,$ i.e. $\r$ has a simple  
zero along each component).
 
\noindent II.\quad $\ov u$ has a pole along $\tilde{F}_1.$
Then each point $\a_i=\tilde{F}_1\cap G_i,$\ $i=1,\dots,q,$
 should be a singular
point of $\ov u,$ and
$\ov{U}_w$ should pass through each $\a_i.$
{}From $(\ov{U}_w,\ov{F}_z)=1$, it
 follows that there is only one component $G_1$ and
the fiber
$\rho^{-1}(0)$ is connected in this case.

Then $S\simeq \BC^2$ (see, for example, \cite{Su}), and is evidently
isomorphic to a
hypersurface.

\noindent III.\quad $\ov u$ is not constant along $\tilde{F}_1.$
Because  $(\ov{U}_w,\tilde{F}_1)=1$ for a general $w,$
it takes every value only once along $\tilde{F}_1.$
{ }From $G_i\cap G_j=\emptyset,$ 
it follows that
 $u_i\ne u_j$ for
$i\ne j,$\
$i,j=1,\dots,s.$

Consider a polynomial $p(u)=(u-u_1)\dots (u-u_q)$
 and $\ov v=\fc{p(\ov u)}{\rho}.$

Since $\ov u$ is finite along $\tilde{F}_1,$
 $\ov v$ is regular and finite at all the points of $S$ and it 
 has a simple
pole along  $\tilde{F}_1.$

Let $A_j=H_j+G_j$ be the divisor $u=u_j.$ Since
 $(\ov{U}_w,\tilde{F}_1)=1$
for a general $w,$
 $(A_j, \tilde {F}_1)=1$ and  $(H_j,\tilde {F}_1)=(A_j, \tilde {F}_1)
-(G_j,\tilde {F}_1)=0.$ Thus, $\tilde {F}_1$ 
does not intersect  zeros of function $\ov v.$ 
 In particular,  the 
intersection points $ s_j=G_i\cap \tilde{F}_1$
are not singular for $\ov v$,  restriction $\ov v \big|_{G_i}$
 has simple poles in $s_j$
and  is linear along each $G_i,$\ $i=1,\dots,q$ (i.e.,
takes every value
$z\in\BP^1$ precisely at one point of $G_i).$

The restriction of $\ov v$ on $S$ we denote by $v,$\ $v\in O(S).$

We define a regular map $\phi: S\to \BC^3$ as
$\phi(s)=(\rho(s),v(s),u(s)),$ and we want to show
that $\phi$ is an isomorphism of $S$ onto a hypersurface
$$S'=\{(x,y,t)\in \BC^3| xy=p(t)\}\subset\BC^3.$$

A)\ $\phi$ is an embedding.
Indeed, the functions $\rho$ and $u$ divide points in $(S\setminus (\cup
G_i)),$ since $\rho$
divides fibers of a line pencil, and $u$ is linear along each fiber
$\rho^{-1}(z),$\ $z\ne 0.$

The values $u\bigm |_{G_i}=u_i$ provide the distinction between the components $G_i$
of $\rho^{-1}(0),$
since $u_i\ne u_j$ when $i\ne j.$

The function $v$ is linear along each $G_i,$ so its values are different in
the different points
of each $G_i.$

B)\ $\phi$ is onto.
Let $s'=S' $ and $s'=(x',y',t').$  If $x'\ne 0,$ then in the fiber 
$\rho^{-1}(x')$ there is a point, 
such that $u(s)=t'.$ (Indeed, $\rho^{-1}(x') \cong\BC$
and $u \bigm |_{\rho^{-1}(x')}$  is linear.) 
Now, $v(s)=\fc{p(u)}{\rho}=\fc{p(t')}{x'}=y',$ so $\phi(s)=s'.$

If $x'=0,$ then $p(t')=0,$ so $t=u_j$ for some $1\le j\le q.$
The function $v$ is linear along the component $G_j,$ so there is a point
$s\in G_j$ such that
$v(s)=y'.$
Then
$\phi(s)=(0,y',u_j)=(0,y',t')=s'\in S'.$\hfill $\qed$ \edm

We proceed with the proof of Theorem 1.
Any surface $S\in \CH$ is a hypersurface by Lemma 5.
If $S\in \CA$ but $S\notin \CH,$ then, by Lemma 4, $[K(S)]\ne0,$
 and, by Lemma 1,
it cannot be
isomorphic to a hypersurface.\hfill $\qed$\edm
 An
example of  a surface $S\in\CA\setminus\CH$ was given
in the Introduction:
 $S\subset \BC^4$  is defined by equations
\smallskip
$\cases xy=(z^2-1)z\\ zu=(y^2-1)y\\ xu=(y^2-1)(z^2-1)\endcases$
\smallskip

 We will show below that this surface is not
 isomorphic to a hypersurface.
 On the other hand, there are two locally 
nilpotent derivations ($lnd$)
defined in the ring
$O(S),$ namely: 
\medskip
$\cases \partial_1x=0\\ \partial_1z=x^2\\ \partial_1y=(3z^2-1)x\\
\partial_1u=2z(y^2-1)x+2y(z^2-1)(3z^2-1)\endcases$
\medskip
$\cases \partial_2u=0\\ \partial_2y=u^2\\ \partial_2z=(3y^2-1)u\\
\partial_2x=2y(z^2-1)u+2z(y^2-1)(3y^2-1)\endcases$
\bigskip
\noindent It follows that $AK(S)=\BC.$

\proclaim{Corollary  of Lemma 1}
The surface $S\subset \BC^4$ defined by equations
\smallskip
$\cases xy=(z^2-1)z\\ zu=(y^2-1)y\\ xu=(y^2-1)(z^2-1)\endcases$
\smallskip
\flushpar is not isomorphic to a hypersurface.
\ep

\demo{Proof} Consider the 2-form $w=\fc{dx\wedge dz}{x}.$
It is regular in  Zariski open subset
 $U_0=\{(x,y,z,u)\in S| x\ne 0\},$ where $(x,z)$ are the
local coordinates.

The fiber $\{x=0\}$ consists of 4 components:
$$\alignat2 & G_1=\{x=0,z=1\};\quad && G_2=\{x=0,z=-1\};\\
&G_3=\{x=0,z=0,y=1\};\quad && G_4=\{x=0,z=0,y=-1\}.\endalignat$$
We consider  Zariski open neighborhoods $U_1,\ U_2,\ U_3,\ U_4$,
respectively, of these components:

$U_1=\{(x,y,z,u)\in S|\ z\ne 0,\ z\ne-1\}$ with local coordinates
$\vp_1=\fc{z-1}{x},$\
$\psi_1=x;$

$U_2=\{(x,y,z,u)\in S|\ z\ne0,\ z\ne 1\}$ with local coordinates
$\vp_2=\fc{z+1}{x},$\
$\psi_2=x;$

$U_3=\{(x,y,z,u)\in S|\ z^2\ne 1,\ y\ne0,\ y\ne-1\}$ with local coordinates
$\vp_3=\fc{y-1}{z},$\ $\psi_3=z;$

$U_4=\{(x,y,z,u)\in S|\ z^2\ne 1,\ y\ne0,\ y\ne 1\}$ with local coordinates
$\vp_4=\fc{y+1}{z}$,\ $\psi_4=z.$

Rewriting $\omega$ in these coordinates, we obtain:
$$\alignat2 &\om=\fc{dx\wedge dz}{x}\quad&&\text{in}\ U_0\\
&\om=d\psi_1\wedge d\vp_1 \quad&&\text{in}\ U_1\\
&\om=d\psi_2\wedge d\vp_2 \quad&&\text{in}\ U_2\\
&\om=-\fc{\psi_3d\vp_3\wedge d\psi_3}{\vp_3\psi_3+1} \quad&&\text{in}\ U_3\\
&\om=-\fc{\psi_4d\vp_4\wedge d\psi_4}{\vp_4\psi_4-1} \quad&&\text{in}\
U_4\endalignat$$

Since $\vp_3\psi_3+1=y\ne0$ in $U_3$ and $\vp_4\psi_4-1=y\ne0$ in $U_4,$
this form is
holomorphic everywhere on $S.$

But $\om\bigm |_{G_3}=\om\bigm |_{G_4}=0$ and 
the divisor $(\om)=G_3+G_4$ is not equivalent
to zero on $S$ by
Lemma 3. So, by Lemma 2, the surface $S$ cannot be isomorphic to a
hypersurface.\hfill$\qed$ \edm
\bk
\subheading{4.\ Corollaries for cylinders and $\BC^+$-actions}

\proclaim{Theorem 2}
Let $S_1,$ $S_2$ be smooth affine surfaces, such that $S_1\in \CH$ and
$S_2\in \CA\setminus \CH.$
Then $S_1\times \BC^k\not\simeq S_2\times \BC^k$ for any $k\in\BN.$\ep

\demo{Proof of Theorem 2}
Assume, to the contrary, that $S_1\times\BC^k\simeq S_2\times \BC^k=W.$
Since $S_1\in \CH,$ by Theorem 1, it is isomorphic to a hypersurface
$S\subset \BC^3,$ and
$W\simeq S\times \BC^k$ is a hypersurface in $\BC^{k+3}$ as well.
By Lemma 1, for any holomorphic $(k+2)$-form $\om$ on $W,$ the divisor
$(\om)$ has to be linearly
equivalent to zero, i.e., there should be a rational
function $f_{\om}$ such
that $(\om)=(f_{\om}).$

Take any holomorphic 2-form $\om'$ on $S_2.$  
Let $z_1,...,z_k$ be  some global coordinates in $\BC^k.$
 Then the form $\om =\om'\wedge dz_1\cdots\wedge dz_k$ is a 
holomorphic  ${k+2}$-form on $W.$ Hence,  $(\om)=(f_{\om})$ 
for a function 
$f_{\om}.$ But then $[(\om')]=[(f_{\om}\bigm|_{S_2})]=0.$ 
It follows that $[K(S_2)]=0$ and, due to Theorem 1, $S_2\in\CH.$
\hfill $\qed$\edm

\proclaim{Theorem 3} A surface $S\in \CA$ admits a fixed-point $\BC^+$-action
if and only if $S\in \CH.$
\ep

\demo{Proof of Theorem 3}
Let $S\in \CA$ and let $\vp_\rho$ be a fixed-point free $\BC^+$-action, let
$\rho$ be a
corresponding line pencil and let $\rho^{-1}(0)$ consist of $q$ components
$G_1,\dots, G_q.$
Consider another surface $S_q=\{xy=(z-1)\dots (z-q)\}\subset \BC^3.$
This  surface is   smooth, affine  
and  has two $\BC^+$-actions; namely,
$$\vp_x^\lm(x,y,z)=(x, \fc{(z+\lm x-1)\dots (z+\lm x-q)}{x}, z+\lm x)$$
and
$$\vp_y^\lm(x,y,z)=(\fc{(z+\lm y-1)\dots (z+\lm y-q)}{y}, y, z+\lm y).$$
Thus, $S_q\in \CA.$
The actions $\vp_x^\lm$ and $\vp_y^\lm$  have no fixed points, 
because the corresponding $lnd'$s,
namely:

$$\partial_x:\partial_x(x)=0,\ \partial _x(z)=x,\ \partial_x(y)=p'(z);$$

$$\partial_y:\partial_y(y)=0,\ \partial_y(z)=y,\ \partial_y(x)=p'(z);$$
never vanish.

The fibers of the line pencil $\rho_x$ 
in $S_q$
corresponding to $\partial_x$ are
the curves
$\{x=\const\}.$
All of them are connected except the fiber $x=0$, which has $q$
connected components.

The fibers of the line pencil $\rho$ in $S$ have
 precisely the same structure.

By the Theorem of Daniliewski and  Fieseler (\cite{D}, \cite {Fi}  )
 the cylinders 
$S\times \BC\simeq
S_q\times \BC.$

But $S_q$ is a hypersurface, hence $S_q\in \CH$ by Theorem 1.
Due to Theorem 2, $S\in \CH$ as well.

 Therefore, if  $S$ admits a fixed-point-free
$\BC^+$-action, then $S\in\CH.$

Now, assume that $S\in\CH.$ As it was shown in Lemma   5, it is 
isomorphic to the surface 
$$S'=\{(x,y,z)\in \BC^3| xy=p(t)\}\subset\BC^3.$$
Since $S$ is smooth, all the roots ${t_1,...t_q}$ of $p(t)$
are simple. 
That is why the $lnd$ $\partial,$  defined as
$$\partial:\partial(x)=0,\ \partial (t)=x,\ \partial(y)=p'(t);$$
does not vanish on $S'.$ 
But then the $\BC^+$-action, defined by $\partial$ has no fixed points.
\hfill $\qed$\edm
\bk

\subheading{Acknowledgments} 
We thank V. Lin and M. Zaidenberg for the idea of  the proof of Lemma 1.
It is our pleasure to thank M. Gizatullin for the discussions concerning
quasihomogeneous  surfaces, and M. Miyanishi and R.V. Gurjar
for the most helpful discussions and examples.

\baselineskip 15pt

\noindent Tatiana M. Bandman, Dept. of Mathematics \& CS,
Bar-Ilan University, Ramat-Gan, 52900, Israel,
 e-mail: bandman\@macs.biu.ac.il.

 \noindent Leonid Makar-Limanov, Dept. of Mathematics \& CS,
 Bar-Ilan University, Ramat-Gan, 52900, Israel, e-mail:
 lml\@macs.biu.ac.il; Dept. of Mathematics, Wayne State University,
 Detroit, MI 48202, USA, e-mail: lml\@math.wayne.edu.

\end
\Refs
\widestnumber \key{KKMLR}


\ref \key {B} \by A. Beauville \book Complex algebraic surfaces
London Math. Soc. Lecture notes \vol 66 \yr 1983 \endref

\ref \key{BML1}\by T. Bandman, L. Makar-Limanov\paper
Cylinders over affine surfaces \jour Japan. Jour. Math.
\vol 26 \yr 2000 \pages 208-217\endref



\ref \key {BML2}\by T. Bandman, L. Makar-Limanov\paper
Affine Surfaces with isomorphic cylinders, preprint,
 e-prints 0001067.\endref




 \ref \key {Ber} \by  J.Bertin \paper Pinceaux de droites
et automorphismes des surfaces affines \jour J. Reine 
und Angew. Math.\vol 341 \yr 1983 \pages 32-53\endref


\ref\key {D} \by W. Danielewski \paper On the cancellation problem and
automorphism groups of affine algebraic varieties,
 preprint,
Warsaw \yr 1989
\endref


\ref \key{Fi} \by K.-H. Fieseler \paper On complex affine surfaces with 
$\BC^+$-action \jour Comment. Math. Helvetici \vol 69 \yr1994
\pages 5-27\endref


\ref \key{G1} \by M.H.Gizatullin \paper Invariants of incomplete 
algebraic surfaces obtained by completions
\jour Math. USSR Izvestiya \vol 5 \yr 1971
\pages 503-516\endref


\ref \key{G2} \by M.H.Gizatullin \paper Quasihomogeneous affine surfaces
\jour Math. USSR Izvestiya \vol 5 \yr 1971
\pages 1057-1081\endref


\ref \key{GD} \by M.H.Gizatullin, V.I.Danilov \paper
 Automorphisms of  affine surfaces.I
\jour Math. USSR Izvestiya \vol 9 \yr 1975
\pages 493-534\endref


\ref \key{Ha} \by Hartshorne \book Algebraic Geometry,
\publ Springer-Verlag, New York, Berlin\yr 1977\endref

 
 \ref \key {KKMLR} \by S. Kaliman, M. Koras, L. Makar-Limanov, P. Russell
\paper $C^*$-actions on $C^3$ are linearizable 
\jour ERA-AMS \vol 3 \yr 1997 \pages 63-71\endref


\ref \key {KML}
\by  S. Kaliman, L. Makar-Limanov\paper On the Russell-Koras
 contractible threefolds \jour
Journ. of the Algebraic Geometry\vol 6(2)\yr 1997
\pages  247--268
\endref


\ref \key {ML1} \by L. Makar-Limanov \paper
 Locally nilpotent derivations, a new ring invariant
and applications, preprint\endref


\ref \key{ML2} \by L. Makar-Limanov \paper Cancellation for curves,
preprint\endref


\ref \key {Mi1} \by M.Miyanishi \book Non-complete algebraic surfaces,
Lecture Notes in Math. \vol 857 \publ Springer-Verlag\yr 1981\endref


\ref \key {MiSu} \by M.Miyanishi, T. Sugie\paper 
Affine surfaces containing cylinderlike open set\jour
J. Math. Univ.Kyoto 
 \vol 20 \yr 1980 \pages 11-42\endref 


\ref \key {S} \by M. Snow\paper Unipotent actions on affine space
\inbook Topological methods in Algebraic Transformation groups
Progress in Math. \vol 80 \yr 1989\pages 165-177
\endref


\ref \key {Su} \by T.Sugie \paper Algebraic 
characterization of the affine plane and the affine 3-space    \inbook Topological methods in Algebraic Transformation groups
Progress in Math. \vol 80 \yr 1989\pages 177-190
\endref



\endRefs
